\documentclass{gtart_h}

\def\ifplaintex{\expandafter\ifx\csname documentclass\endcsname\relax}

\def\gtp{{\mathsurround=0pt\it $\cal G\mskip-2mu$eometry \&\ 
$\cal T\!\!$opology $\cal P\!$ublications}}  

\def\recd{{\small Received:\qua\receiveddate\ifx\reviseddate\relax
\else\qquad Revised:\qua\reviseddate\fi\par}} 

\def\lognumber#1{\def\thelognumber{#1}}
\def\volumenumber#1{\def\thevolumenumber{#1}}
\def\volumeyear#1{\def\thevolumeyear{#1}}
\def\papernumber#1{\def\thepapernumber{#1}}
\def\pagenumbers#1#2{\def\startpage{#1}\def\finishpage{#2}}
\def\published#1{\def\publishdate{#1}}

\def\received#1{\def\receiveddate{#1}}
\def\revised#1{\def\reviseddate{#1}}
\def\accepted#1{\def\accepteddate{#1}}

\def\asciiaddress#1{\def\theasciiaddress{#1}}

\long\def\asciiabstract#1{\long\def\theasciiabstract{#1}}

\let\\\par\let\thelognumber\relax\let\thevolumenumber\relax
\let\thepapernumber\relax\let\thevolumeyear\relax\let\startpage\relax
\let\finishpage\relax\let\publishdate\relax\let\receiveddate\relax
\let\reviseddate\relax\let\accepteddate\relax\let\theasciititle\relax
\let\theasciiauthors\relax\let\theasciiaddress\relax
\let\theasciiabstract\relax

\let\theasciiemail\relax

\ifplaintex
\font\logobig=cmssbx10 scaled 3836
\font\logomed=cmssbx10 scaled 2557
\else
\font\logobig=cmssbx10 scaled 4200
\font\logomed=cmssbx10 scaled 2800
\fi

\long\def\makeagttitle{   
\count0=\startpage
\agt\hfill      
\hbox to 45truept{\vbox to 0pt{\vglue -13truept{\logomed A\kern -.37em{\logobig 
T}\kern -.38em G}\vss}\hss}
\break
{\small Volume \thevolumenumber\ (\thevolumeyear)
\startpage--\finishpage\nl
Published: \publishdate}

\vglue .25truein

{\parskip=0pt\leftskip 0pt plus
1fil\def\\{\par\smallskip}{\Large\bf\thetitle}\par\medskip} \vglue
0.05truein

{\parskip=0pt\leftskip 0pt plus 1fil\def\\{\par}{\sc\theauthors}
\par\medskip}%
 
\vglue 0.03truein 

{\small\leftskip 25truept\rightskip 25truept{\bf Abstract}\stdspace\theabstract

{\bf AMS Classification}\stdspace\theprimaryclass
\ifx\thesecondaryclass\relax\else; \thesecondaryclass\fi\par
{\bf Keywords}\stdspace \thekeywords\par}\vglue 7truept

}   

\ifplaintex
\font\phead=cmsl9 scaled 950
\font\pnum=cmbx10 scaled 913
\font\pfoot=cmsl9 scaled 950
\headline{\vbox to 0pt{\vskip -4.5mm\line{\small\phead\ifnum
\count0=\startpage ISSN 1472-2739 (on-line) 1472-2747 (printed)
\hfill {\pnum\folio}\else\ifodd\count0\def\\{ }%
\ifx\theshorttitle\relax\thetitle\else\theshorttitle\fi\hfill{\pnum\folio}
\else\def\\{ and }{\pnum\folio}\hfill\ifx\theshortauthors\relax\theauthors
\else\theshortauthors\fi\fi\fi}\vss}}
\footline{\vbox to 0pt{\vglue 0mm\line{\small\pfoot\ifnum\count0=\startpage
\copyright\ \gtp\hfill\else
\agt, Volume \thevolumenumber\ (\thevolumeyear)\hfill\fi}\vss}}
\else
\font=cmsl9 scaled 1050
\font\lnum=cmbx10 
\font=cmsl9 scaled 1050
\makeatletter
\def\@oddhead{{\small\lhead\ifnum\count0=\startpage ISSN 1472-2739 
(on-line) 1472-2747 (printed)\hfill {\lnum\number\count0}\else\ifodd\count0
\def\\{ }\ifx\theshorttitle\relax \thetitle \else\theshorttitle\fi\hfill
{\lnum\number\count0}\else\def\\{ and }{\lnum\number\count0}
\hfill\ifx\theshortauthors\relax 
\theauthors\else\theshortauthors\fi\fi\fi}}\def\@evenhead{\@oddhead}
\def\@oddfoot{\small\lfoot\ifnum\count0=\startpage\copyright\ \gtp\hfill\else
\agt, Volume \thevolumenumber\ (\thevolumeyear)\hfill\fi}
\def\@evenfoot{\@oddfoot}
\makeatother
\fi
\let\maketitlepage\makeagttitle
\let\makeshorttitle\maketitlepage
\let\maketitle\maketitlepage

\newwrite\gtoutfile
\long\gdef\makeheadfile{  
{\def\\{, }\def\s{ }
\immediate\openout\gtoutfile head.xxx
\immediate\write\gtoutfile{Proxy-for: \ifx\theasciiauthors\relax
\theauthors\else\theasciiauthors\fi\s<\ifx\theasciiemail\relax\theemail\else\theasciiemail\fi>}
\immediate\write\gtoutfile{\noexpand\\}
\immediate\write\gtoutfile{Authors: \ifx\theasciiauthors\relax
\theauthors\else\theasciiauthors\fi}
{\def\\{ }\immediate\write\gtoutfile{Title: \ifx\theasciititle\relax
\thetitle\else\theasciititle\fi}}
\immediate\write\gtoutfile{Subj-class: GT or SG, GR etc}
\immediate\write\gtoutfile{MSC-class: \theprimaryclass\ifx\thesecondaryclass\relax\else, \thesecondaryclass\fi}
\immediate\write\gtoutfile{Journal-ref: Algebr. Geom. Topol. \thevolumenumber\s
(\thevolumeyear) \startpage-\finishpage}
\immediate\write\gtoutfile{Comments: Published by Algebraic and
Geometric Topology at}
\immediate\write\gtoutfile{\s\s\s  http://www.maths.warwick.ac.uk/agt/AGTVol\thevolumenumber/agt-\thevolumenumber-\thepapernumber.abs.html}
\immediate\write\gtoutfile{\noexpand\\}
\immediate\write\gtoutfile{}
\ifx\theasciiabstract\relax
\immediate\write\gtoutfile{\theabstract}\else
\immediate\write\gtoutfile{\theasciiabstract}\fi
\immediate\write\gtoutfile{}
\immediate\write\gtoutfile{\noexpand\\}
\immediate\write\gtoutfile{}
\immediate\closeout\gtoutfile}}  

\def\maketitlepage{\makeagttitle\makeheadfile}
\let\makeshorttitle\maketitlepage
\let\maketitle\maketitlepage

\lognumber{7}
\volumenumber{5}
\volumeyear{2005}
\papernumber{7}
\pagenumbers{119}{128}
\received{10 December 2004} 
\revised{24 January 2005}
\accepted{27 January 2005}
\published{7 February 2005}

\usepackage[cmtip,arrow]{xy}
\usepackage{pb-diagram,pb-xy}
\usepackage{amsmath,amssymb}

\newtheorem{thm}{Theorem}[section]
\newtheorem*{AT}{Anick's Theorem}
\newtheorem{lem}[thm]{Lemma}
\newtheorem{prop}[thm]{Proposition}
\newtheorem{cor}[thm]{Corollary}
\theoremstyle{definition}

\theoremstyle{remark}
\newtheorem{rmk}[thm]{Remark}
\newtheorem{ex}[thm]{Example}

\newcommand{\Zloc}[1]{\mathbf{Z}_{({#1})}}
\newcommand{\Zmod}[1]{\mathbf{Z}_{{#1}}}

\begin{document}

\title[Hopf algebras up to homotopy and the BSS]{Hopf algebras up to homotopy
and the\\Bockstein spectral sequence}
\author{Jonathan Scott}
\address{Institut de G\'eom\'etrie, Alg\`ebre, et Topologie \\
  \'Ecole Polytechnique F\'ed\'erale de Lausanne \\
  1015 Lausanne, Switzerland}
\asciiaddress{Institut de Geometrie, Algebre, et Topologie\\Ecole 
Polytechnique Federale de Lausanne\\1015 Lausanne, Switzerland}

\email{jonathan.scott@epfl.ch}

\begin{abstract}
Anick proved that every $q$-mild Hopf algebra up to homotopy is
isomorphic to a primitively-generated chain Hopf algebra.  We
provide a new proof, that involves extensive use of the Bockstein
spectral sequence.
\end{abstract}

\asciiabstract{%
Anick proved that every q-mild Hopf algebra up to homotopy is
isomorphic to a primitively-generated chain Hopf algebra.  We
provide a new proof, that involves extensive use of the Bockstein
spectral sequence.}



\maketitle

\section{Introduction}
Let $R$ be a subring of the rationals containing $1/2$. A
\emph{Hopf algebra up to homotopy}, or \emph{Hah}, is an augmented
chain $R$-algebra $A$, equipped with an algebra morphism $\Delta:A
\rightarrow A \otimes A$ that is homotopy-coassociative and
homotopy-cocommutative.  We will suppose that the augmentation
provides a strict co-unit for $\Delta$.  A Hah morphism is an
algebra morphism that commutes with the diagonals up to chain
algebra homotopy.

The aim of this paper is to present an alternate proof of the following
theorem of D. Anick.

\begin{AT}
{\rm\cite[Theorem 5.6]{anick:89}}\qua Let $\rho > 2$ be the least
non-invertible prime in $R$, and let $q \geq 1$. If $A$ is a Hah
whose underlying algebra is free and generated by the range of
dimensions $q$ through $q\rho-1$, inclusive, then $A$ is
isomorphic as a Hah to a primitively-generated chain Hopf algebra.
\end{AT}

Anick's theorem was one of the more important recent developments
in rational homotopy theory. Applied to the Adams-Hilton model of
a space $X$, it allows rational homotopy theorists to apply their
techniques to decidedly non-rational problems.  Indeed, as
outlined in~\cite{halperin:92}, the singular cochain complex
$C^{*}(X;R)$ is then weakly equivalent to a commutative cochain
algebra that plays the role of the algebra of polynomial
differential forms, $A_{PL}(X)$, in rational homotopy theory.

The proof provided by Anick is quite technical, and involves
constructing a ``partial inverse'' to the diagonal map. In the
present paper, the author proposes a proof that uses nothing more
than the Bockstein spectral sequence as detailed in
Browder~\cite{browder:61}, and the results on Hopf algebras from
Milnor and Moore~\cite{milnor-moore:65}.

The key to the new proof of Anick's theorem is
Theorem~\ref{thm:differential-extension} below, where we consider
the result of adjoining a free algebra variable to a
primitively-generated chain Hopf algebra. First we introduce some
terminology. A \emph{free monogenic extension} of a chain algebra
$A$ is a chain algebra $(A \amalg T(x), dx = b)$ where $\amalg$
denotes the free product of algebras, $b \in A_{\deg{x} -1}$ is a
cycle, and $T(x)$ is the tensor algebra on the graded module
$R\{x\}$. If $A$ is a cocommutative chain Hopf algebra, then a
\emph{free monogenic extension of Hopf algebras up to homotopy} is
a free monogenic extension of $A$ along with a
homotopy-coassociative, homotopy-cocommutative choice of reduced
diagonal $\Phi \in (A \otimes A)_{\deg{x}}$ for $x$. A good
reference for the homotopy theory of chain algebras is Anick's
paper~\cite{anick:89}. Effectively, to say that $\Phi$ is a
``homotopy-coassociative'' and ``homotopy-cocommutative'' choice
means that there exist elements $f \in A \otimes A \otimes A$ and
$g \in A \otimes A$ such that $(\Delta \otimes 1 + 1 \otimes
\Delta)\Phi = df$ and $(\tau - 1) \Phi = dg$. If $\Phi=0$ then we
call the extension \emph{trivial}. An isomorphism of extensions
\[
  \theta:(A \amalg T(x), dx=b, \bar{\Delta}x=\Phi)
    \xrightarrow{\cong} (A \amalg T(x), dx=b', \bar{\Delta}x=\Phi')
\]
is a chain algebra morphism that restricts to the identity on $A$,
that satisfies $\theta(x) - x \in A$, and that commutes
with diagonals up to a chain algebra homotopy vanishing on $A$.

Recall that a graded algebra $A$ is called $q$-\emph{reduced} if
$A_{n}=0$ for $0 < n < q$.

\begin{thm}\label{thm:differential-extension}
Let $(A,\partial)$ be a connected, primitively generated,
$q$-reduced, $R$-free chain Hopf algebra of finite type. Suppose
$A \rightarrow A \amalg T(x)$ is a free extension of Hopf algebras
up to homotopy, with $\partial x$ primitive.  If $\deg x < q\rho$,
then the extension is isomorphic to the trivial one.
\end{thm}

Theorem~\ref{thm:differential-extension} will be proved in
Section~\ref{sec:proof}.

\begin{proof}[Proof of Anick's Theorem]
Let $A = TV$ be a Hopf algebra up to homotopy.
Let $\{ v_{j} \}$ be a well-ordered basis of $V$, where we have chosen an
ordering such that $\deg{v_{i}} < \deg{v_{j}}$ implies that $i < j$.
Suppose inductively that $A_{(i-1)} = T(v_{1}, \ldots v_{i-1})$ is a
primitively generated Hopf algebra.  Adding the generator $v_{i}$
yields the chain algebra extension
\[
  A_{(i-1)} \rightarrow A_{(i-1)} \amalg T(v_{i})
\]
where $\bar{\Delta}v_{i}=\Phi$.  It is easy to see that $[\partial
v_{i}]$ is primitive in $H(A_{(i-1)})$.  By
Lemma~\ref{lem:HL-in-HUL} below, $H_{m}(P(A_{(i-1)} \otimes
\Zmod{p})) \cong PH_{m}(A_{(i-1)} \otimes \Zmod{p})$ if $m <
q\rho$, for all primes $p$ in $R$.  Since $P(A_{(i-1)} \otimes
\Zmod{p}) = P(A_{(i-1)}) \otimes \Zmod{p}$ in these degrees, we
conclude that $H_{m}(PA_{(i-1)}) \cong PH_{m}(A_{(i-1)})$, and so
we may choose $\partial v_{i}$ to be primitive in $P(A_{(i-1)})$.
Now we may apply Theorem~\ref{thm:differential-extension}.
\end{proof}

\begin{rmk}
Anick assumes a homotopy counit, but then goes on to show that every
Hah is isomorphic to one with a strict
counit~\cite[Lemma 5.4]{anick:89}.
\end{rmk}

We use the usual definitions for chain algebras, coalgebras, Hopf
algebras, and homotopies of chain algebra morphisms~\cite[Section
2]{anick:89}. Algebras are assumed to be connected and augmented
to the ground ring.  The linear dual of a chain complex $V$ is
denoted $V^{*}$.

The author would like to take this opportunity to thank
Jean-Claude Thomas and Jim Stasheff, for useful comments and
suggestions, as well as the referee, whose recommendations led to
a major reworking of the article.  The end result is more
elementary, less technical, more accessible, and just plain
better.

\section{The Proof of
Theorem~\ref{thm:differential-extension}}\label{sec:proof}

Let $p$ be an odd prime, and denote by $\Zmod{p}$ the field of
integers modulo $p$. We prove a proposition about cocommutative
extensions of Hopf algebras over $\Zmod{p}$.  We then prove
several lemmas that allow us to manage the Bockstein spectral
sequence, before proving the inductive step,
Proposition~\ref{prop:induction}.  We finish the section with the
proof of Theorem~\ref{thm:differential-extension} proper.

We begin by quickly reviewing some facts about Hopf algebras,
gleaned from the standard reference of Milnor and
Moore~\cite{milnor-moore:65}.  Let $C$ be a connected,
associative, coassociative, graded Hopf algebra with
multiplication $\mu:C \otimes C \rightarrow C$, diagonal $\Delta:C
\rightarrow C \otimes C$ and augmentation $\varepsilon:C
\rightarrow R$. Let $I(C) = \ker\varepsilon$. The \emph{space of
indecomposables} of $C$ is the quotient $I(C)/\mu(I(C)\otimes
I(C))$.  Let $\bar{\Delta}:I(C) \rightarrow I(C) \otimes I(C)$ be
the reduced diagonal, defined by $\bar{\Delta}c = \Delta c - c
\otimes 1 - 1 \otimes c$.  Then $P(C) = \ker{\bar{\Delta}}$ is the
\emph{subspace of primitive elements}.  The natural map $I(C)
\rightarrow Q(C)$ restricts to define a map $P(C) \rightarrow
Q(C)$, natural with respect to morphisms of Hopf algebras.

Denote by $TC(V)$ the \emph{tensor coalgebra} on $V$.  As a graded
module, $TC(V) = T(V) = \oplus_{n \geq 0}V^{\otimes n}$.  Elements
of $V^{\otimes n}$ are traditionally denoted
$[v_{1}|\cdots|v_{n}]$.  The diagonal is defined by
\[
    \Delta[v_{1}|\cdots|v_{n}]
        = \sum_{j=0}^{n}[v_{1}|\cdots|v_{j}]\otimes
        [v_{j+1}|\cdots|v_{n}].
\]
In particular, $P(TC(V))=V$.  Note that, if $V$ is finite type and
free as a graded $R$-module, then $(TV)^{*} = TC(V^{*})$ as a
graded coalgebra.

While the following proposition likely qualifies as folklore, we
provide a proof nonetheless.

\begin{prop}\label{prop:extension-mod-p}
Let $A$ be a connected, primitively-generated Hopf algebra of
finite type over $\Zmod{p}$. Suppose that
\[
  A \rightarrow A \amalg T(x)
\]
is a coassociative and cocommutative extension. If $A$ is
$(r-1)$-connected and $\deg(x) < rp$, then $A \amalg T(x)$ is
primitively generated.
\end{prop}

\begin{proof}
Set $B=A \amalg T(x)$.  Dualize the extension $A \rightarrow B$ to
obtain a morphism of connected commutative Hopf algebras $f:B^{*}
\rightarrow A^{*}$. It suffices to show that the natural morphism
$P(B^{*}) \rightarrow Q(B^{*})$ is injective.  Let $\xi:B^{*}
\rightarrow B^{*}$ be the Frobenius map defined by $\xi(x) =
x^{p}$. Then $\xi B^{*}$ is a sub Hopf algebra of $B^{*}$.
By~\cite[Proposition 4.21]{milnor-moore:65}, there is an exact
sequence
\[
    0 \rightarrow P(\xi B^{*}) \rightarrow P(B^{*})
    \rightarrow Q(B^{*})
\]
so it suffices to show that $P(\xi B^{*})=0$.

Let $V = Q(A)$ and $W = Q(B)$, so that $W = V \oplus
\Zmod{p}\{x\}$. Let $\sigma:V \rightarrow A$ be a splitting of the
natural projection $\pi:I(A) \rightarrow V$.  Extend $\sigma$ to a
splitting $\tau:W \rightarrow B$ of $\rho:I(B) \rightarrow W$.
Then $\sigma$ and $\tau$ extend to algebra epimorphisms $TV
\rightarrow A$ and $TW \rightarrow B$ such that the diagram of
algebra morphisms
\[
\begin{diagram}
    \node{TV} \arrow{e,V} \arrow{s,l,A}{\sigma} \node{TW}
        \arrow{s,r,A}{\tau} \\
    \node{A} \arrow{e,V} \node{B}
\end{diagram}
\]
commutes. Dualize to obtain the commutative diagram of coalgebra
morphisms,
\[
\begin{diagram}
    \node{B^{*}} \arrow{e,t,A}{f} \arrow{s,l}{\tau^{*}}
        \node{A^{*}} \arrow{s,r}{\sigma^{*}} \\
    \node{TC(W^{*})} \arrow{e,t,A}{g} \node{TC(V^{*})}
\end{diagram}
\]
wherein the vertical arrows are monomorphisms of coalgebras. Then
the morphisms $P(\sigma^{*}):P(A^{*}) \rightarrow P(TC(V^{*})) =
V^{*}$ and $P(\tau^{*}):P(B^{*}) \rightarrow P(TC(W^{*}))$ $= W^{*}$
are also monic. Furthermore, $P(g):W^{*} \rightarrow V^{*}$ is the
canonical projection that kills $x^{*}$, the basis element dual to
$x \in W$. Suppose $b^{p} \in P(\xi B^{*}) \subset P(B^{*})$. Then
$P(\tau^{*})(b^{p}) = \alpha x^{*} + v^{*}$, according to the
direct sum decomposition $W^{*} = \Zmod{p}\{x^{*}\} \oplus V^{*}$.
Then $0 = P(\sigma^{*})P(f)(b^{p}) = P(g)P(\tau^{*})(b^{p}) =
P(g)(\alpha x^{*} + \beta v^{*}) = v^{*}$, so $v^{*}=0$. Therefore
$P(\tau^{*})(b^{p})=\alpha x^{*}$.  But $\deg b^{p} \geq rp > \deg
x^{*}$ so $\alpha =0$. Since $P(\tau^{*})$ is a monomorphism,
$b^{p}=0$ and the result follows.
\end{proof}

Let $A$ be a chain Hopf algebra over $\Zmod{p}$; that is, $A$ is
simultaneously a Hopf algebra and a chain complex, such that the
differential $\partial$ is both a derivation and a coderivation.
Then $H(A)$ is naturally a Hopf algebra. Let $i_{A}:P(A)
\rightarrow A$ be the natural inclusion of the primitive subspace.
Then $H(i_{A}):H(PA) \rightarrow H(A)$ factors through the
inclusion $i_{H(A)}:PH(A) \rightarrow H(A)$ to yield a morphism
$j:H(PA) \rightarrow PH(A)$.

Let $A$ be a chain Hopf algebra such that $\partial = 0$ on
$A_{<n}$.  If $a \in A_{n}$ is not a cycle, then $a$ is
indecomposable, and $\partial a$ is primitive.

Let $B$ be a connected, primitively-generated Hopf algebra of
finite type over $\Zmod{p}$.  Since $p$th powers vanish in
$B^{*}$, they also vanish in the commutative Hopf algebra
$H(B^{*}) = H(B)^{*}$.  By~\cite[Proposition
4.20]{milnor-moore:65}, it follows that $H(B)$ is primitively
generated.

\begin{lem}\label{lem:HL-in-HUL}
Let $A$ be a connected, primitively-generated chain Hopf algebra
of finite type over $\Zmod{p}$. If $A$ is $q$-reduced, then
$j:H(PA) \rightarrow PH(A)$ is an isomorphism in degrees $< qp$.
\end{lem}

\begin{proof}
The primitive filtration $\{ F_{k} \}$ on $A$ leads to a
first-quadrant spectral sequence of commutative,
primitively-generated Hopf algebras that converges to $E^{0}H(A)$.
By~\cite[Proposition 5.11]{milnor-moore:65}, $p$th powers vanish
outside of bidegree $(0,0)$, $P(E^{0}A)=E^{0}_{1,t}A=P(A)_{1+t}$,
and $P(E^{0}A) \xrightarrow{\cong} Q(E^{0}A)$. The morphism $H(PA)
\rightarrow PH(A)$ corresponds to the edge homomorphism
$E^{1}_{1,t}A \rightarrow E^{\infty}_{1,t}A$.

Since $P(E^{0}A)$ is a sub chain complex of $E^{0}A$, we may write
$P(E^{0}A)=X \oplus Y \oplus Z$, with $d^{0}:X \xrightarrow{\cong}
Y$ and $H(P(E^{0}A))\cong Z$.  It follows that $E^{0}A$ is the
tensor product of differential Hopf algebras of the following
forms:
\begin{enumerate}
    \item $B_{1}=(\wedge(z),0)$,
    \item $B_{2}=(\Zmod{p}[z]/(z^{p}), 0)$,
    \item $B_{3}=(\wedge(x) \otimes \Zmod{p}[y]/(y^{p}),dx=y)$,
    \item $B_{4}=(\Zmod{p}[x]/(x^{p}) \otimes \wedge(y),dx=y)$,
\end{enumerate}
where $\wedge(-)$ denotes the exterior algebra.  The
indecomposables all lie in $E^{0}_{1,*}A$. A calculation shows
that $H(B_{3})=\wedge(xy^{p-1})$ and $H(B_{4})=\wedge(x^{p-1}y)$.
It follows that $Q(E^{1}A)$ is concentrated in columns $1$ and
$p$, and $P(E^{1}A) \xrightarrow{\cong} Q(E^{1}A)$.  Let $d^{r}$
be the first non-vanishing differential in the spectral sequence,
$r \geq 1$. Let $u$ be an element of lowest total degree such that
$d^{r}u \neq 0$.  Then $u \in Q(E^{r}A)$ and $d^{r}u \in
P(E^{r}A)$. Thus $u \in E^{r}_{p,t}A$ for some $t$, $d^{r}u \in
E^{r}_{1,t+p-2}A$, and so $r=p-1$.  An element of $F_{p}$ with
non-vanishing differential has total degree at least $qp+1$, so $t
\geq qp + 1 - p$. Therefore $d^{p-1}u \in E^{p-1}_{1,qp - 1}A$, so
the edge homomorphism $E^{1}_{1,t}A \rightarrow E^{\infty}_{1,t}A$
is an isomorphism for $t < qp - 1$.  We conclude that $H(P(A))
\rightarrow P(H(A))$ is an isomorphism in degrees $<qp$.
\end{proof}

\begin{ex}
Consider the commutative, primitively-generated Hopf algebra over
$\Zmod{p}$, $A = \Zmod{p}[x] \otimes \wedge(y)$, with $\deg{x} =
2n$, and $n \geq 1$. The differential is the derivation determined
by the rule $\partial x = y$.  Then $PA = \Zmod{p}\{y, x, x^{p},
x^{p^{2}}, \ldots \}$. It follows that $H(PA)=\Zmod{p}\{ x^{p},
x^{p^{2}}, \ldots \}$. On the other hand, $H(A) = \Zmod{p}[x^{p}]
\otimes \wedge(x^{p-1}y)$, so $PH(A) = \Zmod{p}\{ x^{p-1}y, x^{p},
x^{p^{2}}, \ldots \}$. So $j:H(PA) \rightarrow PH(A)$ fails to be
surjective in degree $2np-1$.
\end{ex}

\begin{ex}
Let $B = \wedge(x) \otimes \Zmod{p}[y]$, with $\deg{x} = 2n+1$.
The differential is given by $\partial x = y$.  Then $PB =
\Zmod{p}\{ x, y, y^{p}, y^{p^{2}}, \ldots\}$, and so $H(PB)=
\Zmod{p}\{ y^{p}, y^{p^{2}}, \ldots \}$.  Meanwhile, $H(B) = 0$,
so $PH(B)=0$.  Therefore $j:H(PB) \rightarrow PH(B)$ is not
injective.
\end{ex}

Let $A$ be a connected primitively-generated Hopf algebra of
finite type over $\Zloc{p}$, the ring of integers localized at
$p$. Recall from~\cite{browder:61} that the \emph{Bockstein
spectral sequence} of $A$ in homology modulo $p$, denoted
$(E^{r},\beta^{r})$, is a spectral sequence of
primitively-generated Hopf algebras over $\Zmod{p}$ that converges
\[
  E^{1} = H(A;\Zmod{p}) \Rightarrow
    \left(\frac{H(A)}{\mathrm{torsion}}\right) \otimes \Zmod{p}.
\]
We now introduce some notation to avoid confusion when following
elements up and down the terms of the spectral sequence. Let $C$
be a chain complex, and denote by $E^{r}$ its Bockstein spectral
sequence modulo $p$. If $x \in C$ is a cycle modulo $p$ that
survives to the $r\,$th term of the spectral sequence, then we
denote its equivalence class in $E^{r}$ by $[x]_{r}$.  In
particular, $[x]_{2} = [[x]_{1}]$, and so on.

The following elementary lemma  clarifies what it means for one
element to be the higher Bockstein of another.

\begin{lem}\label{lem:bss-rep}
Let $C$ be a chain complex and denote by $E^{r}$ the $r\,$th term
of the Bockstein spectral sequence of $C$ modulo $p$. Suppose
$a_{r}, b_{r} \in E^{r}$ satisfy $\beta^{r}a_{r} = b_{r}$.
Suppose further that $a,b\in C$ represent $a_{r}$ and $b_{r}$,
respectively.  That is, $a$ and $b$ are cycles modulo $p$ that
survive to $E^{r}$, where $[a]_{r}=a_{r}$ and $[b]_{r}=b_{r}$.
Then there exist $c,e \in C$ such that $d(a+pc)=p^{r}(b + pe)$.
\end{lem}

\begin{proof}
We proceed by induction on $r$.  For $r=1$, if $a$ is a cycle
modulo $p$, then $da = px$ for some $x \in C$.  By definition,
$\beta^{1}[a]_{1} = [x]_{1}$, so $[x]_{1}=[b]_{1}$.  Thus there
exist elements $c,e \in C$ such that $x = b + pe - dc$.  Therefore
$da = px = pb + p^{2}e - d(pc)$ and the statement holds, grounding
the induction.

Assume that the lemma is true for $r = k-1$ for some $k \geq 2$,
and suppose that $\beta^{k}[a]_{k} = [b]_{k}$.  Since $a$ survives
to the $k\,$th term of the spectral sequence,
$\beta^{k-1}[a]_{k-1}=0$.  Applying the inductive hypothesis, we
find that there exist elements $f,g \in C$ such that $d(a + pf) =
p^{k}g$.  By definition, $\beta^{k}[a]_{k}=[g]_{k}$.  Thus
$[g]_{k} = [b]_{k}$, so $[g-b]_{k}=0$, that is, there exists some
$y \in C$ that survives to $E^{k-1}$, where $[g-b]_{k-1} =
\beta^{k-1}[y]_{k-1}$.  By the inductive hypothesis, there exist
elements $z,e \in C$ such that $d(y + pz) = p^{k-1}(g-b-pe)$.
Therefore $d(a+pf) = p^{k}g = p^{k}b + p^{k+1}e + d(py + p^{2}z)$.
Rearranging, we find that $d(a + p(f - y - pz))=p^{k}(b+pe)$.  Set
$c=f-y-pz$ to complete the inductive step and the proof.
\end{proof}

\begin{cor}\label{cor:bss-rep}
If $[b']_{r} = [b']_{r}$ in $E^{r}$, then there exist elements
$e,f \in C$ such that $p^{r-1}b' = p^{r-1}b'' + p^{r}e + df$.
\end{cor}

\begin{proof}
Since $[b']_{r} = [b'']_{r}$, there exists $a \in C$ such that
$\beta^{r-1}[a]_{r-1} = [b'-b'']_{r-1}$.  Apply
Lemma~\ref{lem:bss-rep}.
\end{proof}

\begin{lem}\label{lem:key}
Let $A$ be a connected, $q$-reduced, primitively generated chain
Hopf algebra over $\Zloc{p}$. Suppose that $\partial a = p^{r}b$,
for some $a,b \in A$, and that $\partial w = p^{r-1}\bar{\Delta}b$
for some $w \in A \otimes A$. If $\deg b < qp$, then there exist
elements $x,y \in A$ such that $x \otimes 1 \in P(A \otimes
\Zmod{p})$ and $\partial(a - x - py) = 0$ .
\end{lem}

\begin{proof}
We claim that for each $i \geq 0$, there exist elements $b_{i},
y_{i}, z_{i} \in A$ and $\Psi_{i} \in A \otimes A$, with $z_{i}$
primitive mod $p$, that satisfy
\begin{equation}\label{eq:a}
    \partial(a - z_{i} - py_{i}) = p^{r+i}b_{i}
\end{equation}
and
\begin{equation}\label{eq:w}
    \partial(w - \bar{\Delta}y_{i} - \Psi_{i}) = p^{r+i-1}\bar{\Delta}b_{i}.
\end{equation}
When $i=0$, we take $z_{0}=y_{0}=0$, $b_{0}=b$, and $\Psi_{0}=0$.

Suppose for some $i \geq 0$ that we have the elements $b_{i}$,
$y_{i}$, and $z_{i}$ as above. If $b_{i}$ does not survive to the
$(r+i)\,$th term, then there exist $v, b_{i+1} \in A$ such that
$\partial v = p^{s}b_{i} - p^{s+1}b_{i+1}$ for some $s<r+i$.  Then
$\partial(a - z_{i} - py_{i} - p^{r+i-s}v)=p^{r+i+1}b_{i+1}$ and
$\partial(w-\bar{\Delta}(y_{i} + p^{r-s-1}v) - \Psi_{i}) =
p^{r+i}\bar{\Delta}b_{i+1}$. Set $z_{i+1} = 0$, $y_{i+1} = y_{i} +
p^{r+i-s}v$ and $\Psi_{i+1}=\Psi_{i}$.  We note that $\deg b_{i+1}
= \deg b_{i}$.

Suppose on the other hand that $b_{i}$ survives to $E^{r+i}$. Then
$[b_{i}]_{r+i}$ is primitive by equation (\ref{eq:w}). Let $\{
[b]_{r+i} \}$ denote the homology class of $[b]_{r+i}$ in
$H(PE^{r+i})$. Since $\beta^{r+i}[a - z_{i}]_{r+i} =
[b_{i}]_{r+i}$, it follows that $\{[b]_{r+i}\} \in
\ker\{j:HPE^{r+i} \rightarrow PE^{r+i+1}\}$. By
Lemma~\ref{lem:HL-in-HUL}, there exists $[z]_{r+i} \in PE^{r}$
such that $[b]_{r+i} = \beta^{r+i}[z]_{r+i}$.  By repeated
application of Lemma~\ref{lem:HL-in-HUL}, we may suppose that $z$
is primitive modulo $p$;  say, $\bar{\Delta}z=p\Psi$.  By
Lemma~\ref{lem:bss-rep},
\begin{equation}\label{eqn:thing}
  \partial(z + py) = p^{r+i}(b_{i} - pb_{i+1})
\end{equation}
for some elements $y,b_{i+1} \in A$.  Substituting (\ref{eq:a}) in
(\ref{eqn:thing}) and rearranging, we find that
\[
  \partial\left(a - (z_{i}+z) - p(y_{i} + y)\right) = p^{r+i+1}b_{i+1}.
\]
Taking the reduced diagonal of each side of (\ref{eqn:thing}) and
using (\ref{eq:w}), we obtain
\[
  \partial\left(w - \bar{\Delta}(y_{i}+y) - (\Psi_{i}+\Psi)\right) =
    p^{r+i}\bar{\Delta}b_{i+1}.
\]
Set $y_{i+1} = y_{i} + y$, $z_{i+1}=z_{i}+z$ and $\Psi_{i+1}+\Psi$
to complete the induction.

Now, since $H_{\deg b}(A)$ is finitely generated, its torsion
submodule has an exponent, say $m$, so that the maximum order of
torsion is $p^{m}$.  In particular, there exists a chain $u \in A$
such that $\partial u = p^{m}b_{m-r+1}$.  By construction,
$\partial(a - z_{m-r+1} - py_{m-r+1}) = p^{m+1}b_{m-r+1}$.
Therefore $\partial(a - z_{m-r+1} - p(y_{m-r+1} + u))=0$,
completing the proof.
\end{proof}

\begin{prop}\label{prop:induction}
Let $A$ be a connected, primitively generated, $q$-reduced,
$\Zloc{p}$-free chain Hopf algebra of finite type. Let $n < qp$.
If $\Phi \in (A \otimes A)_{n}$, $f \in (A \otimes A \otimes
A)_{n+1}$, and $g \in (A \otimes A)_{n+1}$ satisfy $\partial \Phi
= 0$, $\partial f = (\bar{\Delta} \otimes 1 - 1 \otimes
\bar{\Delta}) \Phi$, and $\partial g = (\tau - 1) \Phi$, then for
all $r \geq 1$, there exist cycles $a_{r} \in A_{n}$, $\Phi_{r}
\in (A \otimes A)_{n}$, and a chain $\Omega_{r} \in (A \otimes
A)_{n+1}$, such that
\[
  \bar{\Delta}a_{r} = \Phi - p^{r}\Phi_{r} + \partial\Omega_{r}.
\]
\end{prop}

\begin{proof}
We proceed by induction on $r$.  For $r=1$, we note that $\Phi \in
H(A;\Zmod{p}) \cong E^{1}$ is strictly coassociative and
cocommutative.  Since $\deg\Phi < qp$,
Proposition~\ref{prop:extension-mod-p} tells us that there exists
$[a'_{1}]_{1} \in E^{1}$ such that $\bar{\Delta}[a'_{1}]_{1} =
[\Phi]_{1}$. Thus there exist chains $\Psi'_{1}, \Omega_{1} \in A
\otimes A$ such that
\begin{equation}\label{eq:a1}
  \bar{\Delta}a'_{1} = \Phi - p\Phi'_{1} + \partial\Omega_{1}.
\end{equation}
Since $a'_{1}$ is a cycle modulo $p$, $\partial a'_{1} = pb_{1}$
for some $b_{1} \in A$. Differentiating (\ref{eq:a1}), we find
that $\bar{\Delta}b_{1}=-\partial\Phi'_{1}$. Now apply
Lemma~\ref{lem:key}, to obtain elements $x \in P(A \otimes
\Zmod{p})$, $y \in A$, such that $\partial(a'_{1} - x - py) = 0$.
Set $a_{1} = a'_{1} - x - py$. Since $x$ is primitive modulo $p$,
$\bar{\Delta}x = p\Psi$ for some $\Psi$. Then $\bar{\Delta}a_{1} =
\Phi - p(\Phi'_{1} + \Psi + \bar{\Delta}y) + \partial \Omega_{1}$.
Set $\Phi_{1} = \Phi'_{1} + \Psi + \bar{\Delta}y$.  Since
$\partial\Phi_{1} = -\bar{\Delta}\partial a_{1} = 0$, the
induction is grounded.

Suppose now that
\begin{equation}\label{eq:rubba}
  \bar{\Delta}a_{r-1} = \Phi - p^{r-1}\Phi_{r-1} + \partial\Omega_{r-1}.
\end{equation}
Apply $\bar{\Delta} \otimes 1 - 1 \otimes \bar{\Delta}$ to
(\ref{eq:rubba}) and rearrange to obtain
\[
  \partial(f +
  (\bar{\Delta} \otimes 1 - 1 \otimes \bar{\Delta})\Omega_{r-1})
  = p^{r-1}(\bar{\Delta} \otimes 1 - 1 \otimes \bar{\Delta})
  \Phi_{r-1},
\]
so
\[
  \beta^{r-1}[f]_{r-1} = (\bar{\Delta} \otimes 1 - 1 \otimes \bar{\Delta})
    [\Phi_{r-1}]_{r-1}.
\]
Similarly, by applying $\tau - 1$ to (\ref{eq:rubba}) and
rearranging, we obtain
\[
  \beta^{r-1}[g]_{r-1} = (\tau - 1) [\Phi_{r-1}]_{r-1},
\]
so $[\Phi_{r-1}]_{r}$ is a candidate for a diagonal in $E^{r}$. By
Proposition~\ref{prop:extension-mod-p}, there exists
$\tilde{a}_{r} \in E^{r}$ such that
$\bar{\Delta}\tilde{a}_{r}=[\Phi_{r-1}]_{r}$. There exists a
representative $\tilde{a} \in A$ of $\tilde{a}_{r}$ and an element
$\tilde{b} \in A$, such that $\partial\tilde{a}=p^{r}\tilde{b}$.
By Corollary~\ref{cor:bss-rep},
\begin{equation}\label{eq:zubba}
  p^{r-1}\bar{\Delta}\tilde{a}
  = p^{r-1}\Phi_{r-1} - p^{r}\Phi'_{r} + \partial{u}
\end{equation}
for some $\Phi'_{r},u \in A \otimes A$. Adding equations
(\ref{eq:rubba}) and (\ref{eq:zubba}), we obtain
\[
  \bar{\Delta}(a_{r-1}+p^{r-1}\tilde{a})
  = \Phi - p^{r}\Phi'_{r} + \partial(u + \Omega_{r-1}).
\]
Taking boundaries, we find that $-\partial \Phi_{r} =
p^{r-1}\bar{\Delta}\tilde{b}$. So by Lemma~\ref{lem:key}, there
exist elements $z, y \in A$, with $z$ primitive modulo $p$, such
that $\partial(\tilde{a} - z - py)=0$. Set $a_{r} = a_{r-1} +
p^{r-1}(\tilde{a} - z - py)$.  Write $\bar{\Delta}z = p\Psi$. Then
\[
    \bar{\Delta}a_{r} = \Phi - p^{r}(\Phi'_{r} + \Psi +
    \bar{\Delta}y) + \partial(u + \Omega_{r-1}).
\]
Set $\Phi_{r} = \Phi'_{r} + \Psi + \bar{\Delta}y$ and $\Omega_{r}
= u + \Omega_{r-1}$.  Then $\partial\Phi_{r} = -
\bar{\Delta}\partial a_{r} = 0$, so the inductive step is
complete.
\end{proof}

\begin{proof}[Proof of Theorem~\ref{thm:differential-extension}]
Consider $C = \Omega A / \Omega^{\geq 3} A$, where $\Omega(-)$ is
the Adams cobar construction~\cite{adams:56}. Since $A$ is finite
type, $H_{n-2}(C)$ is a finitely generated $R$-module. We note
that $R$ is a principal ideal domain, so there is a finite list of
primes $p$ for which $H_{n-2}(C)$ contains $p$-torsion.  For each
of those primes, we make the following argument.

Let $\Phi = \bar{\Delta}x$, and let $f$ and $g$ be the relevant
homotopies.  Since $\partial x$ is primitive, it follows that
$\partial\Phi=0$.  Apply Proposition~\ref{prop:induction}. Then
$[\Phi] \in H_{n-2}(C)$ and $[\Phi] = p^{r}[\Phi_{r}]$ for all
$r$. Since $H_{n-2}(C)$ is finitely generated, it follows that
$[\Phi]=0$. Therefore there exist $\Psi \in A \otimes A$ and $a
\in A$ such that $\bar{\Delta}a = \Phi + \partial \Psi$. Define
\[
  \theta : (A \amalg T(x), \bar{\Delta}x = \Phi)
    \rightarrow (A \amalg T(x), \bar{\Delta}x = 0)
\]
by $\theta|_{A} = 1_{A}$ and $\theta(x)=x+a$.  Then $\theta$ is an isomorphism
of chain algebras that commutes with the diagonals up to homotopy.
\end{proof}

\Addresses\recd

\end{document}